\documentclass[a4paper]{article}

\usepackage[english]{babel}
\usepackage[T1]{fontenc}
\usepackage{amsmath}
\usepackage{amssymb}
\usepackage[dvips]{graphicx}
\usepackage{subfigure}
\usepackage{cases}
\usepackage{empheq}
\usepackage{bbm}
\usepackage{color}
\usepackage{psfrag}
\usepackage{mathtools,verbatim,url}

\begin{document}
\title{Elastic wave propagation in complex geometries: A qualitative comparison between two high order finite difference methods}
\author{Kristoffer Virta, Christopher Juhlin, Gunilla Kreiss}
\date{}
\maketitle

\begin{abstract}
We compare two high order finite-difference methods that solve the elastic wave equation in two dimensional domains with curved boundaries and material discontinuities. Two numerical experiments are designed with focus on wave boundary interaction, the response of a pressure wave impinging on a circular cavity and the wave field generated by an explosive impact on the wall an underground tunnel. Qualitative comparisons of the results are made where similarities and differences are pointed out.  
\end{abstract}

\section{Introduction}
Numerical methods that solve the elastic wave equation are important tools for the understanding of elastic wave propagation. In particular when studying the effect of geometry and material properties in a domain of interest as geometrical features and inhomogeneous materials in general makes the construction of exact solutions impossible. The response of a signal of known frequency and duration can serve as a diagnostic tool to understand effects of layering in the Earth as well as predicting the location and composition of underground mineral deposits and reservoirs of gas and oil in the Earth's interior. Numerical methods can also be used to simulate the effect of seismic events such as earthquakes.

Finite difference (FD) modeling for wave propagation problems is a widely used numerical approach. In the bulk of the computational domain away from boundaries FD methods are easily constructed to high order accuracy and resolution requirements are in general known in terms of grid points per wavelength \cite{Kreiss:1972}. The computational grid used with FD methods grid is structured and thus enables highly efficient parallel implementations \cite{Bohlen:2002, SW4, Henshaw:2016}. Boundaries of the domain introduces boundary conditions and geometrical features to the domain. Both boundary conditions and geometries require numerical treatment that guarantees both accuracy and stability. In the context of the FD method one way to achieve geometric flexibility is to solve the elastic wave equation in curvilinear coordinates \cite{Petersson:2015,Duru_Kreiss_Mattsson}. The advantage of  this approach is that geometries are now represented highly accurate by the computational grid but the computations become more involved as the transformation to curvilinear coordinates introduces variable coefficients and in the curvilinear coordinate system the geometry is represented by anisotropic wave motion. Another approach is to ignore conditions at interfaces and curved boundaries and proceed in the original Cartesian coordinate system. The drawback of this approach, albeit less technical, is that geometries are not represented as accurately and the correct boundary and interface conditions can not be guarantied to be satisfied to order of accuracy. In fact, these conditions are only satisfied to first order accuracy \cite{Brown:1984}.                

FD modeling has been used in a variety of geophysical applications, Bohlen et.al \cite{Bohlen:ch5} used numerical modeling to study the effects of shape and composition of large ore bodies on the scattered wave field. In the context of determining the feasibility of storage of nuclear waste an understanding of the influence of sub surface structures such as layering and fracturing on interior motions was gained thorough FD simulations in \cite{Gritto:2004}. A distinct feature of elastic wave propagation is the complex wave boundary interaction phenomena that can occur, at a traction free boundary surface waves appears. Numerical simulations with FD methods that involve surface waves and their conversion to interior waves has been used to predict geological structures ahead of tunnels \cite{Jetschny:2011,Luth:2008}.

In this paper we make a qualitative comparison between two FD methods that solves the elastic wave equation but treats geometries and boundary conditions differently. One of the methods represent geometries by a transformation to curvilinear coordinates and imposes boundary and interface conditions explicitly so that they are satisfied to order of accuracy. The other method does not impose boundary conditions explicitly and represents geometries by material heterogeneities. The focus of the present study is on wave boundary interaction. In a first problem we compare the effect of an incoming pressure wave on the boundary of a circular cavity and in a second problem we compare the wave fields resulting from an explosive impact on the wall of an underground tunnel. 
\section{Numerical modeling}
The purpose of the present study is to compare a well established numerical method for geophysical problems to one more recently developed method. Both methods use high order (higher that 2) FD stencils to approximate spatial derivatives and solve the elastic wave equation in two spatial dimensions. The main differences of the methods are how curved geometries and boundary and interface conditions are approximated. We now present details of each numerical modeling scheme.  
\subsection{Standard velocity-stress finite difference modeling}
A commonly used method to solve wave propagation problems in heterogeneous elastic media is the velocity-stress formulation of the elastic wave equation. In this formulation the equation is set up as a first order system. The velocity-stress formulation was used in a numerical solver by Virieux \cite{Virieux:1986} using 2nd order operators in time and space and later by Levander (1988) using 4th order operators in space and 2nd order in time. In this study we use a code developed by Juhlin \cite{Juhlin:1995} that is based on the formulation of Levander \cite{Levander:1988} and available in the Seismic Unix software package \cite{seismic_Unix}. We refer to this formulation as the standard velocity-stress (S-VS)  method since it is widely used. In this formulation a regularly spaced Cartesian grid is used and grid points are required throughout the media being modeled. This implies that cavities within the grid must be assigned velocities and densities representing air otherwise the scheme becomes unstable. As a consequence when geometries such as cavities or material interfaces are present they are represented by discontinuous material properties. At the discontinuity no explicit treatment is given to the arising interface conditions and the geometry can only be resolved to first order. For this reason we do not expect high order accuracy in the vicinity of geometries and material interfaces. For an example grid with cavities and material heterogeneity see the left of Figure \ref{fig:grids}.    
\subsection{SBP-SAT finite difference modeling}
A more recent method is based on the 2nd order displacement formulation of the elastic wave equation and uses Summation-By-Parts (SBP) FD operators \cite{Strand:1994} to approximate spatial derivatives. Boundary conditions such as traction free, radiation and Dirichlet boundary condition as well as internal interface conditions are explicitly imposed to high order of accuracy with the Simultaneous-Approximation-Term (SAT) technique \cite{Carpenter:1994}. The resulting SBP-SAT method is a consistent and provably stable fifth order in space and fourth order in time method that solves the elastic wave equation in two spatial dimensions \cite{Duru_Virta, Duru_Kreiss_Mattsson}.

At an interface between two materials the continuity of both displacements and tractions is imposed. When the boundary of the domain is adjacent to air as in the case of a cavity the boundary is assumed to be free of traction. That is, the domain is assumed to be adjacent to a vacuum. This assumption is considered valid since the density of surrounding material is in general several magnitudes larger than that of air. Geometries of the domain are treated by solving the elastic wave equation in curvilinear coordinates. In this way geometries are represented to high order of accuracy. For an example of a grid representing a domain with cavities and a material heterogeneity see the right of Figure \ref{fig:grids} in the figure it is clear that grid points are aligned with the geometrical features. Due to high order approximations of both boundary conditions and geometries we expect no loss of formal order of accuracy due to the presence of boundaries and interfaces. For a detailed study of the performance of the SBP-SAT method in curved geometries see \cite{Virta:2015}.     

\subsubsection{Numerical verification of the SBP-SAT method}
\label{ss:sb_sat_num_exp}
In the following experiment the geometry is identical to that of the study of Section \ref{s:cylinder_comp} but initial data and outer boundary conditions are chosen such that an analytic solution is known. The purpose of this experiment is to illustrate the correctness of the SBP-SAT method and investigate how to choose the grid size using the SBP-SAT method in the study of Section \ref{s:cylinder_comp}.  

Consider a time-harmonic plane pressure wave ($P$-wave) impinging on a circular cavity of radius $r>0$ centered at the origin. As the wave strikes the boundary of the cavity mode conversion occurs and both $P$-waves and shear waves ($S$-waves) are present in the resulting solution. An analytic solution to this problem is described in \cite{Virta:2015} and software to generate initial data for a numerical experiment is available from \cite{Virta:2015a}. 

\begin{figure}[tb]
	\centering
	\subfigure{
	\psfrag{eee}[][][5][0]{TIME}
	\includegraphics[scale=0.36]{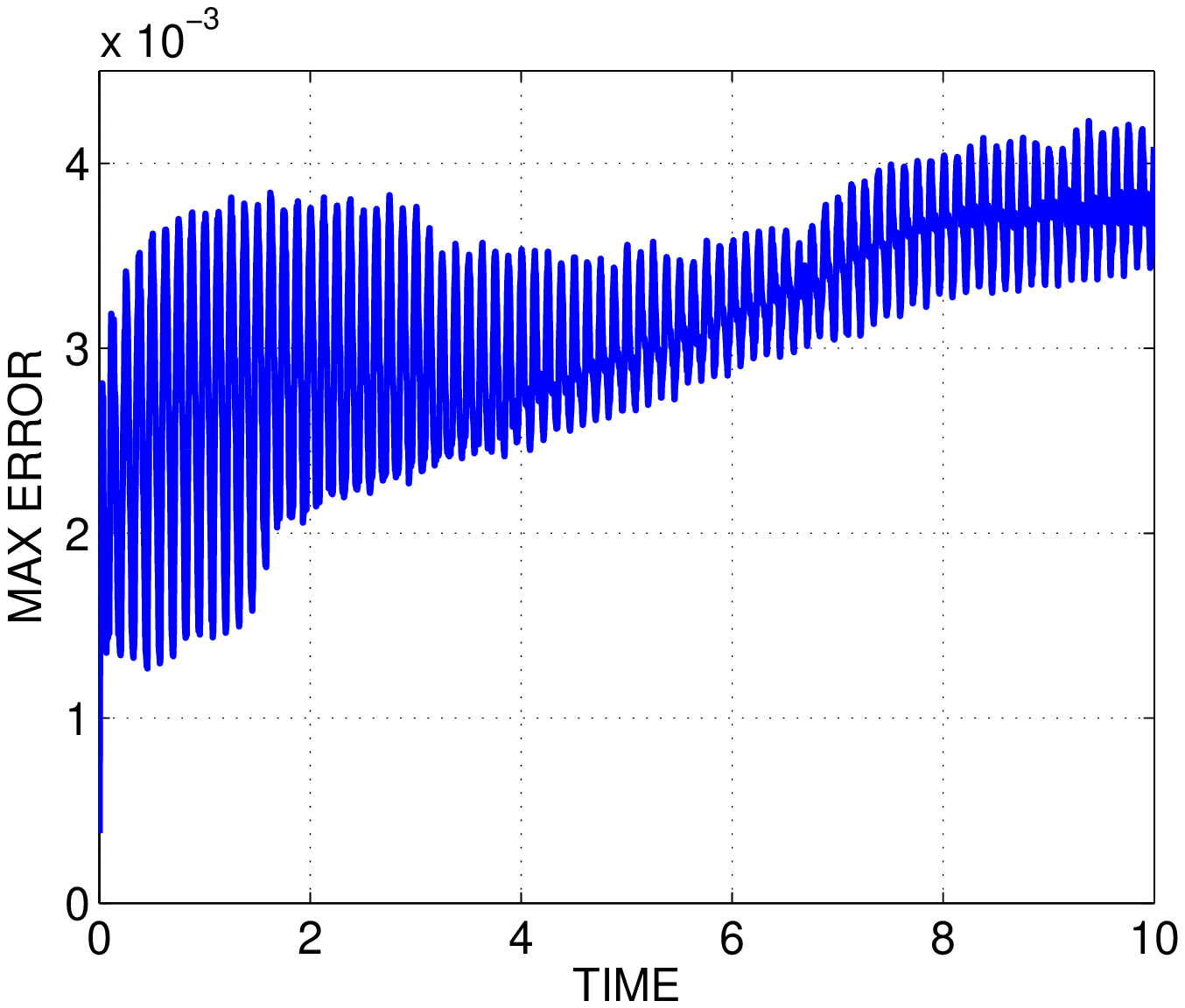}
	}
	\subfigure{
	\includegraphics[scale=0.36]{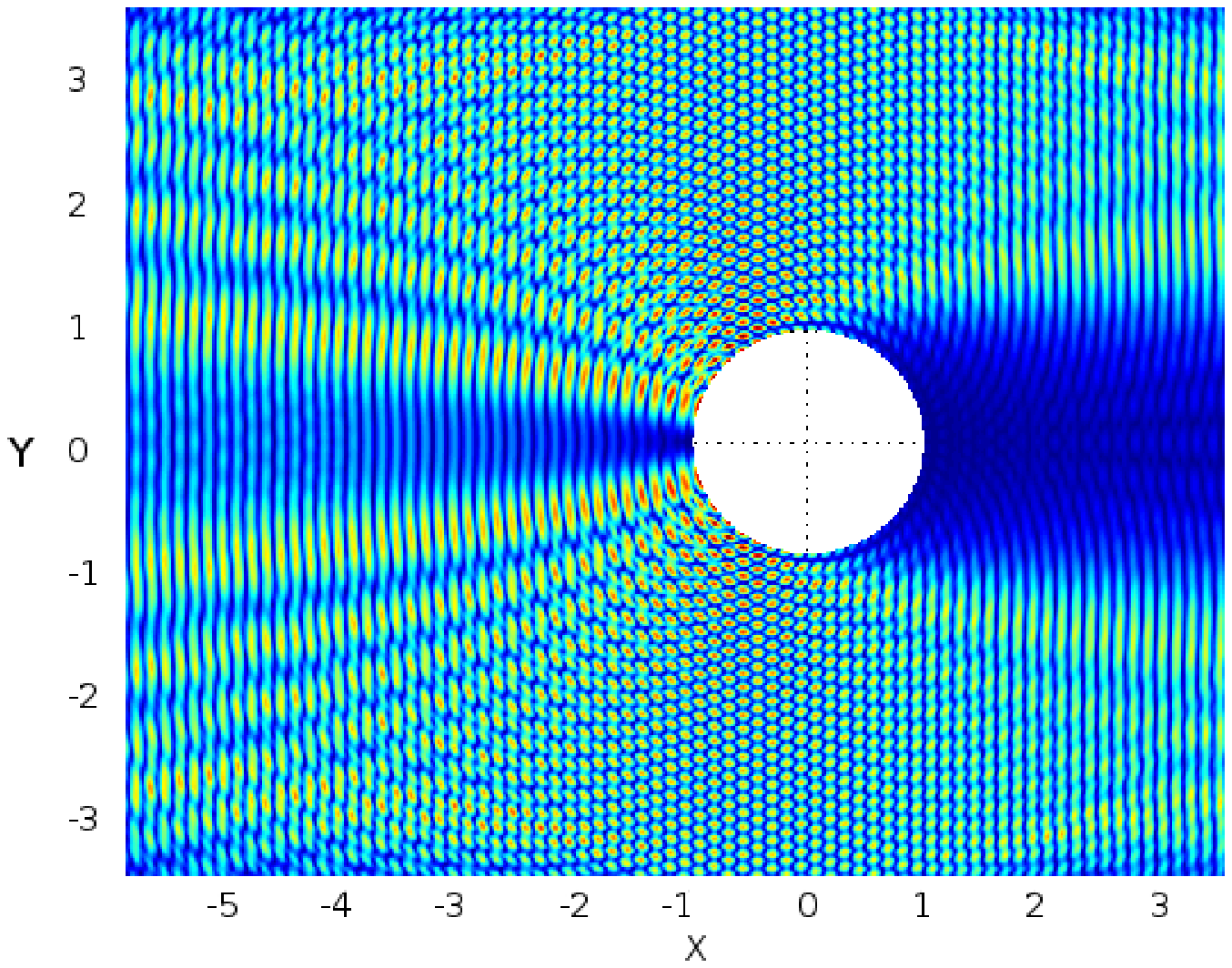}
	}
	\caption{Left: relative maximum error in the numerical solution obtained with the SBP-SAT finite difference method as a function of time. Right: magnitude of the displacement field used for initial data in the experiment of Section \ref{ss:sb_sat_num_exp}.}	
	\label{fig:er_and_sol}
\end{figure}

In the experiment we fix the radius of the circular cavity to be $r=1$. We choose material parameters such that the wavelength of the incoming and scattered $P$-waves is $1/4$ and such that the wavelength of the scattered $S$-waves is $1/8$. The temporal period of the solution is taken to be $1/4$. The right of Figure \ref{fig:er_and_sol} displays the magnitude of the displacement field of the analytic solution used as initial data. We perform a numerical experiment by propagating the solution numerically in the time interval $t \in [0,10]$. The computational domain is $[-5.9,3.6] \times [-3.9,3.9]$ with the circular cavity centered at the origin and is discretized using about $12$ grid points per wavelength of the present $S$-waves, which have the smallest wavelength. The relative maximum error is computed at each time step and plotted in the left of Figure \ref{fig:er_and_sol} where it is seen that a relative maximum error of about $0.5\%$ is achieved with a resolution of $12$ grid points per smallest wavelength. 
         
\section{Seismic waves around a circular cavity}
\label{s:cylinder_comp}
\begin{figure}[tb]
	\centering
	\includegraphics[scale=0.3]{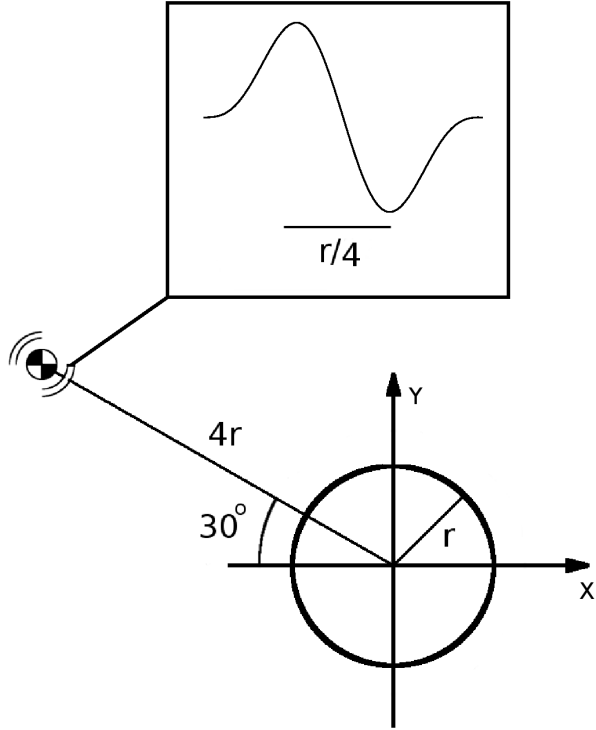}
	\caption{Illustration of the experiment of Section \ref{s:cylinder_comp}}	
	\label{fig:setup}
\end{figure} 
We consider a circular cavity centered at the origin and surrounded by a homogeneous elastic material. The radius of the circular cavity is $r>0$.  At a distance of $4r$ and at an angle of $30^\circ$ from the negative $x$-axis we place a transient source generating $P$-waves of wavelength approximately $r/4$. Horizontal and vertical displacements are recorded at the points $N=(0,r)$, $S=(0,-r)$, $E=(r,0)$ and $W=(-r, 0)$ on the boundary of the cavity in the time interval $t \in [0,10]$. The setup is displayed in Figure \ref{fig:setup}.

The $P$-wave source is obtained by choosing an internal forcing $F$ acting on the vertical and horizontal displacements when using the SBP-SAT method. When the S-VS method is used the vertical and horizontal velocities are subjected to forces. The horizontal and vertical forces are given by the first and second components, respectively, of   
\begin{equation*}
	F(x,y) = \nabla \partial(x-x_s,y-y_s) h(t).
\end{equation*}
%
Here $\partial(x,y)$ is the two dimensional Dirac delta function, 
\begin{equation}
	x_s = 4 r \cos(5/6 \pi), y_s = 4 r \sin(5/6 \pi),
\end{equation}
and 
\begin{equation}	
	h(t) = \left\{
		\begin{array}{ll}
			\sin(2\pi \omega t) - 1/2 \sin(4 \pi \omega t), &t \in (0, 1/\omega)\\
			0, &t  \notin (0,1/\omega)		 
		\end{array},\omega = 2.\right.
\end{equation}
A description of the discretization of the Dirac delta function is given in \cite{Walden:1999}. Initial data is identically zero. We choose material parameters such that $c_p = r$ and $c_s = r/2$, where $c_p$ and $c_s$ are the wave speeds or $P$ and $S$ waves, respectively. This choice gives a Poisson ratio $\nu = 1/3$ which is representative of Aluminum. The amplitude of the Fourier transform of $h(t)$ has local maxima in the vicinity of the points $\pm \omega$ and $\pm 2 \omega$ and decays quadratically for larger frequencies. For this reason we estimate the highest significant frequency in $h(t)$ to be $2\omega$, then the shortest present wavelength of the $P$-waves becomes approximately $c_p/(2 \omega) = r/4$. Initially only $P$-waves are present but when these waves impinge on the boundary of the circular cavity mode conversion occurs and $S$-waves are generated. Due to the choice of $c_s = r/2$ the wavelength of the converted $S$-waves is approximately $r/8$.

The SBP-SAT method models the interior of the circular cavity by a vacuum instead of air. This modeling is considered valid since the density of air inside the  circular cavity is about $2 \times 10^3$ times smaller than the surrounding Aluminum. At the boundary of the cavity a traction free boundary condition  is explicitly imposed. Based on the numerical experiment of Section \ref{ss:sb_sat_num_exp}, the computational domain used with the SBP-SAT method is discretized with at least $12$ grid points per wavelength of the present $S$-waves. The SBP-SAT method uses a grid with grid points located exactly at the boundary of the cavity. We expect that with this resolution the problem is solved with a relative maximum error of about $0.5\%$.

For the modeling of the problem with the S-VS finite difference method a Cartesian grid with grid spacing corresponding to $12$ grid points per shortest wavelength of the waves outside of the circular cavity. Since the velocity and density cannot be zero in the cavity values of $c_p = 0.05 r$ and  $c_s = 0 $ in the cavity, this choice represents an air filled cavity. This implies that waves will propagate through the cavity, rather than only around it. Furthermore, the low velocity in the cavity results in that the wave field is under-sampled within it and that numerical dispersion will influence the solution. At the boundary of the air filled cavity the material properties are discontinuous. In the S-VS method the position of the boundary is modeled only to first order, we expect this to cause a first order error.

Both methods truncate the computational domain with radiation boundary conditions. 
\subsection{Results}
\begin{figure}[htb]
		\centering
		\subfigure{
		\includegraphics[scale=0.75]{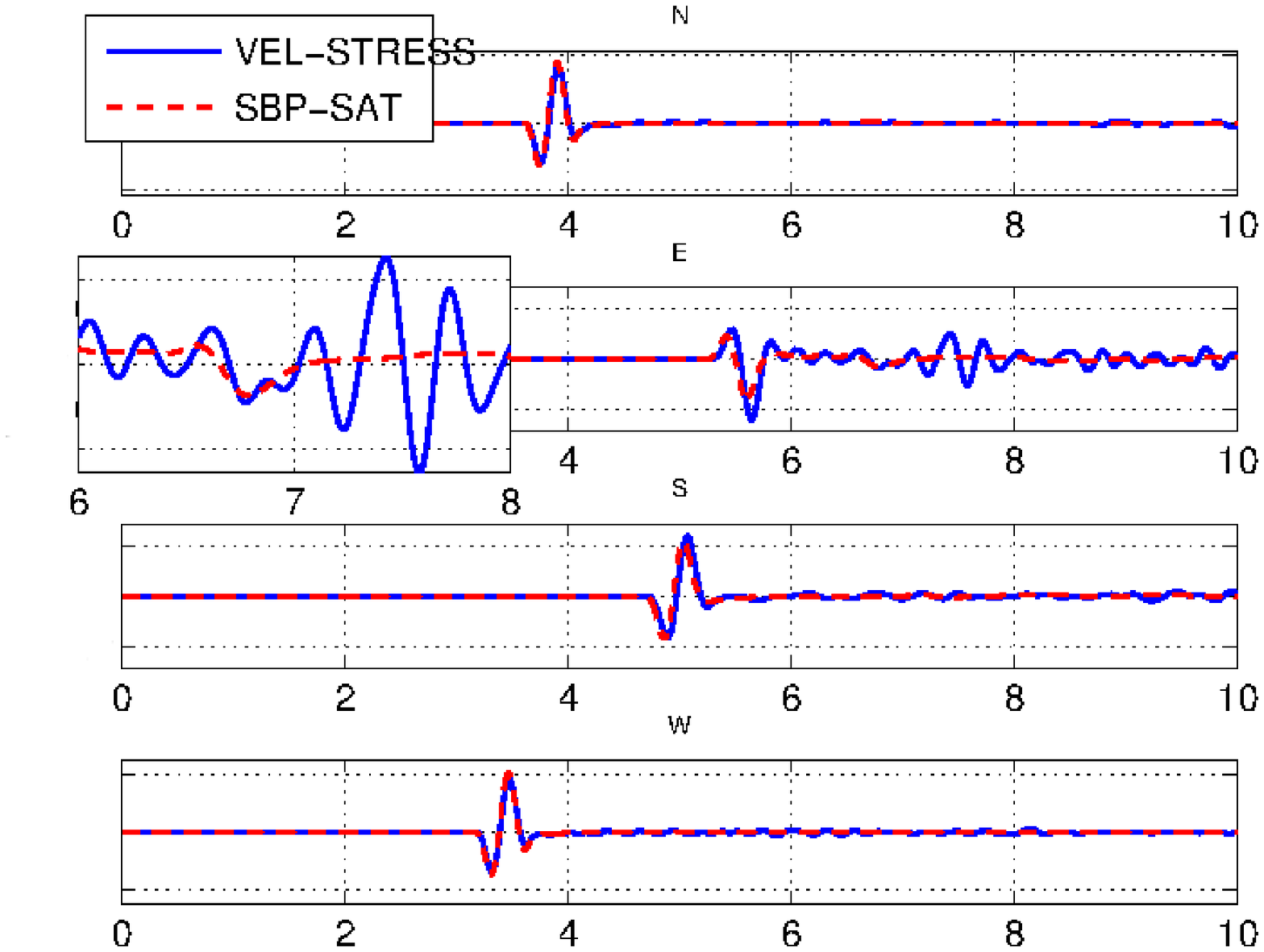}
		}
		\centering
		\subfigure{
		\includegraphics[scale=0.75]{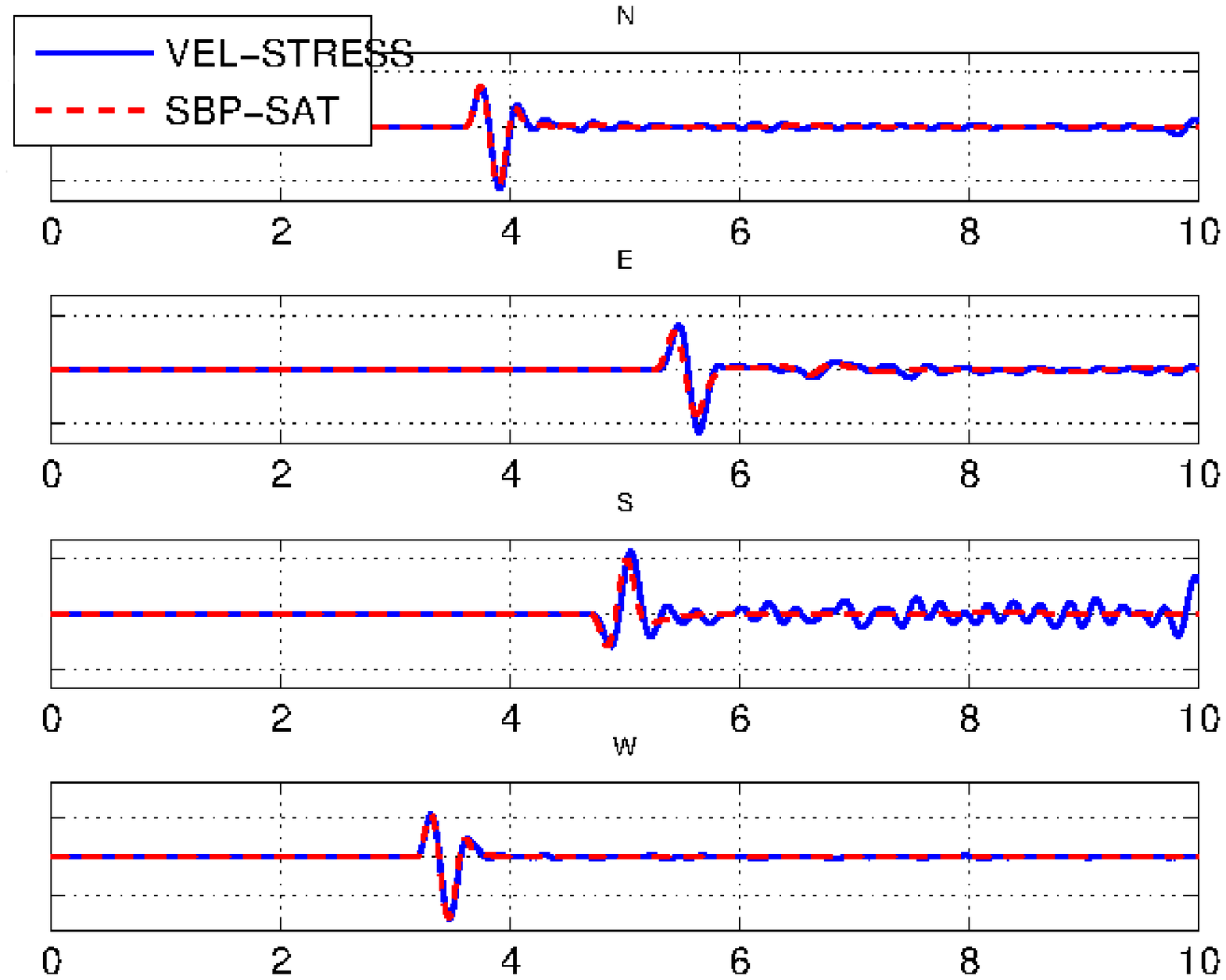}
		}
		\caption{Top: horizontal displacement components. Bottom: vertical displacement components.}
		\label{fig:recordings}
\end{figure}

Recordings of horizontal and vertical displacements at the points $N$, $S$, $E$ and
$W$ on the boundary of the cavity obtained with both numerical methods are plotted as functions of time in Figure \ref{fig:recordings}. At all points the first arrivals are recorded at the same time for the two methods, but at later times differences in amplitudes and arrival times are apparent. The points $N$ and $W$ lie closest to the incoming wave front, and in the recordings at these points the amplitudes and periods of the first arrivals are identical in "eye-norm" for both methods. After the first arrival the S-VS method produces a whole sequence of small waves whereas with
the SBP-SAT method no such waves are not present. When the wave front reaches the points E and S it has already encountered the boundary of the cavity. With the S-VS method waves are traveling both through the air of the cavity and around it, while with the SBP-SAT method no wave
motion is modeled inside the cavity. This difference is seen in figure \ref{fig:recordings} as a difference in amplitude of the first arrival (which corresponds to the wave traveling around the cavity), and significant latter arrivals with the S-VS method. These latter arrivals result from waves traveling through the air of the cavity. The accuracy the latter arrivals is questionable, since the grid size is chosen to resolve the waves in the surrounding medium and wavelengths inside the cavity are much smaller than outside due to a much smaller wave speed there. We expect dispersion errors to affect the waves traveling through the cavity. Other reasons to question the correctness of the latter arrivals is that the jump in material properties at the boundary of the cavity is not handled accurately.

A close up of the recorded horizontal displacement in the time interval $t\in[6,8]$ at the point $E$ is shown as an inset in figure \ref{fig:recordings}. In this time interval a latter arrival corresponding to a converted shear wave is visible in the recording obtained with the SBP-SAT method. This arrival is obscured by waves coming from the interior of the cavity in the recording obtained with the S-VS method. Since the obscuring waves are expected to suffer from dispersion errors we argue that information reveling boundary phenomena may be lost if boundary conditions and geometries are not explicitly treated and all present waves are sufficiently resolved.          

\section{Seismic waves in a general 2D domain}
\begin{figure}[htb]
	\centering
	\subfigure{\includegraphics[scale=0.5]{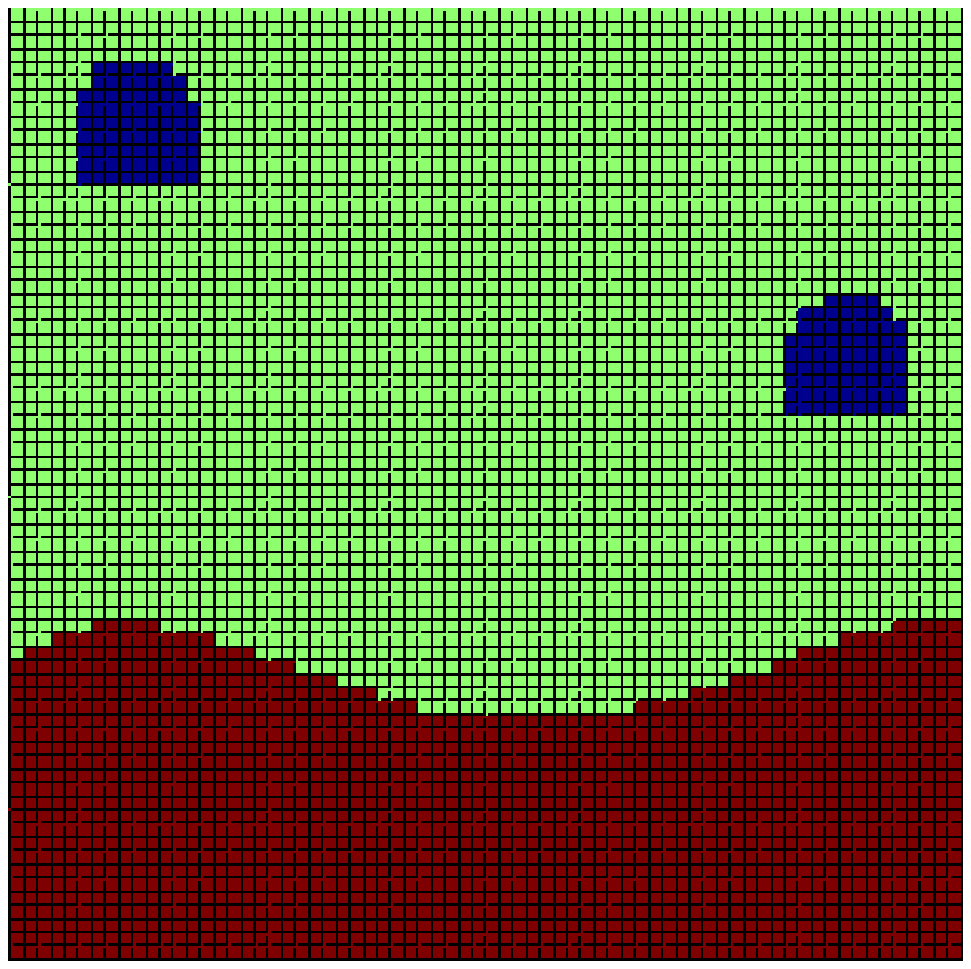}}
	\subfigure{\includegraphics[scale=0.5]{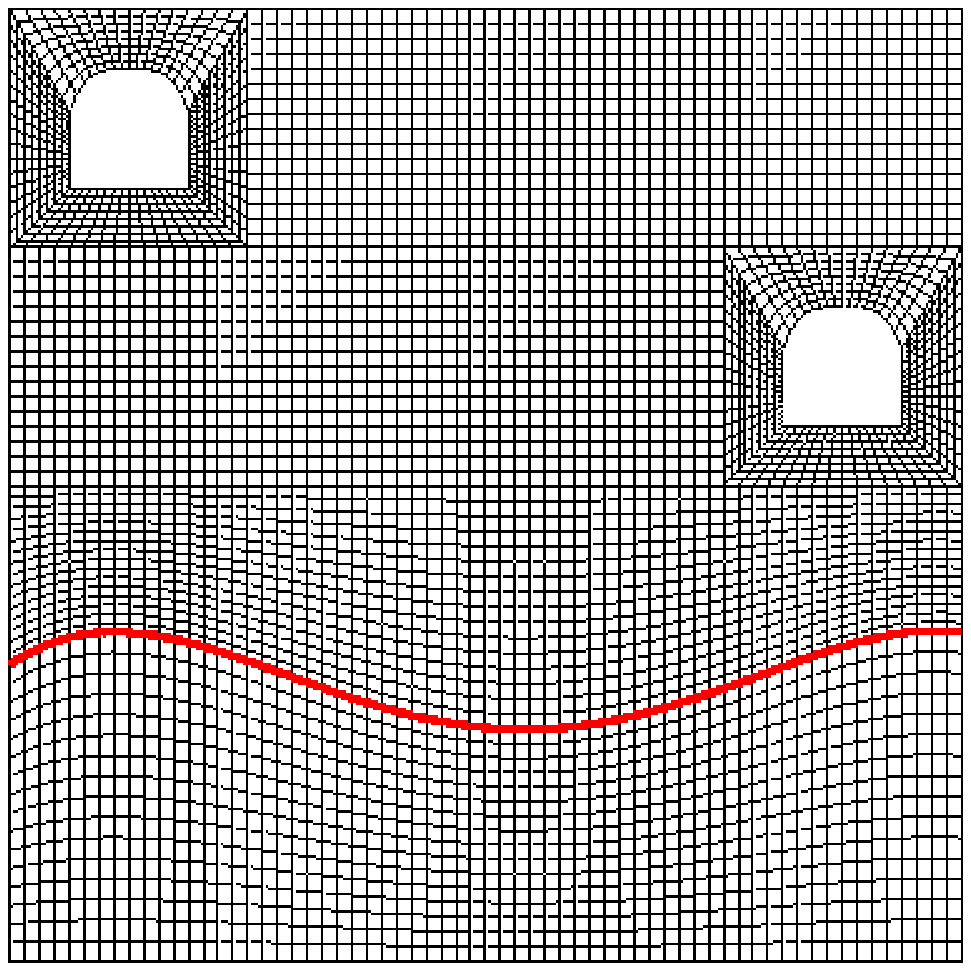}}	
	\caption{Left: Cartesian grid used with the SVS-FD method, different colors indicate different materials. Right: curvilinear grid used with the SBP-SAT FD method, the red curve indicates the interface between the two differing materials.}
	\label{fig:grids}
\end{figure}

As an application we consider a fictional cross section of the Earths interior. The cross section contains two cavities representing tunnels and the medium surrounding the tunnels has material properties with a discontinuity at a curved interface. The geometry of the cross section is displayed in Figure \ref{fig:grids}, where also the grids for the S-VS and SBP-SAT FD methods are visualized. We scale time such that in the upper material surrounding the tunnels we have $c_p = 1.765 m/s$, $c_s =1m/s$ and in the lower material we have $c_p = 1.892 m/s$, $c_s = 1.088m/s$. The Poisson ratio of the upper and lower materials are $0.2636$ and $0.2530$, respectively. This represents typical rocks  in the crust of the earth. The distance between the midpoints of the tunnels is $6.324 m$. At the right wall of the left tunnel waves are generated by an explosive impact normal to the tunnel wall. The purpose of this experiment is to illustrate the geometric possibilities of the used methods and to compare the resulting wave fields visually. The grid used with the SBP-SAT method has grid points aligned with the boundaries of the tunnels and the material interface where boundary and interface conditions are treated explicitly, whereas the S-VS method uses a Cartesian grid with no explicit treatment of boundary and interface conditions. In the SBP-SAT model the interior of the tunnels are considered to be a vacuum due to the large difference in density between air and the surrounding rock. With the S-VS method the interior of the tunnels are assigned material properties representing air.   

\begin{figure}[htb]
	\centering
	\includegraphics[scale=0.25]{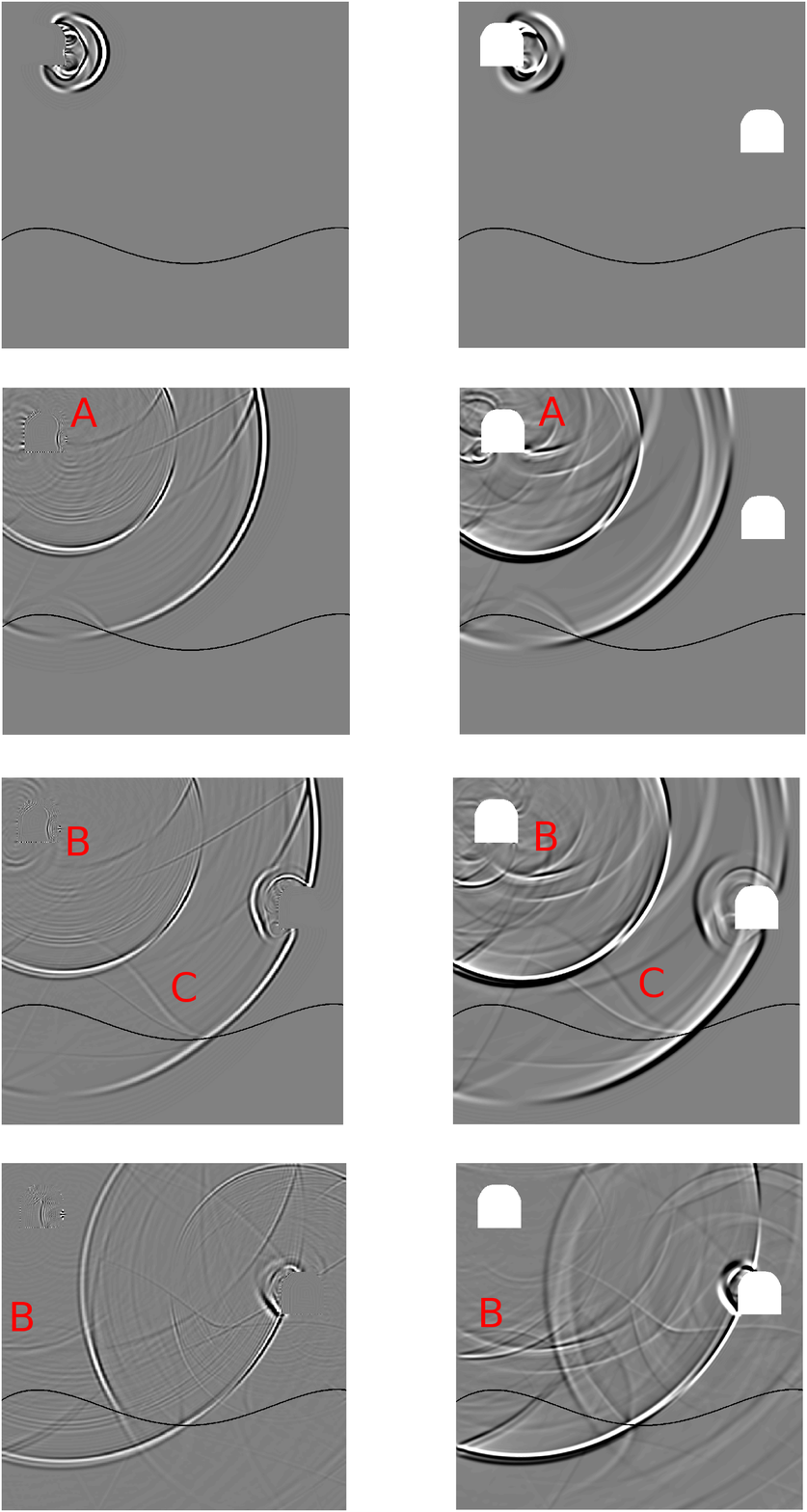}
	\caption{Snapshots of the horizontal component of the displacement field at the times $t=0.6$, $t=2.8$, $t=3.6$ and $t=5.6$. Left: standard velocity-stress finite difference modeling. Right: SBP-SAT finite difference modeling.}
	\label{fig:snaps_x}
\end{figure}

Figure \ref{fig:snaps_x} shows snapshots of the horizontal component of the modeled displacement fields at the times $t=0.6$, $t=2.8$, $t=3.6$ and $t=5.6$ after triggering the explosive impact. The snapshot sequence reveals the location of the direct wave fronts at different times, the locations agree well between the two methods. Rayleigh surface waves are generated at the boundary of the left tunnel, marked by \color{red} A \color{black} in the figure. The Rayleigh waves are more noticeable in the wave field computed with the SBP-SAT method. The surface waves are attenuated due to the curved boundary of the tunnel \cite{Epstein:1976} and parts of the surface waves are being converted into interior waves, marked by \color{red}B \color{black} in the figure. When the wave front impinges on the material interface it is being both reflected and refracted, marked by \color{red} C \color{black} in the figure. The shape of the reflections is affected by the shape of the material interface and agree between both methods. 

We conclude that the numerical treatment of the boundary of the tunnel cavity highly influences the energy content in the present Rayleigh surface waves. As a consequence the waves converted from Rayleigh waves are also less in energy content in the wave field generated by the S-VS method.           
\section{Conclusions}
We have made a qualitative comparison between the results of using two different numerical methods to solve two problems of elastic wave propagation in a domain with non-trivial geometrical features and boundary conditions. The first problem involved the response of an impinging pressure wave on the boundary a circular cavity. The second problem considered the wave field resulting from an explosive impact on the wall of a tunnel in a heterogeneous medium. The used methods were an established widely used standard method based on the first order velocity-stress formulation of the elastic wave equation and a more recent SBP-SAT method  based on the second order displacement formulation of the elastic wave equation. The main difference between the two methods was the treatment of geometries of boundaries and material interfaces and the associated boundary and interface conditions. For the standard method no explicit treatment was done while the more recent method approximates boundary and interface conditions to order of accuracy. Also, the standard method was assigning values of material parameters inside cavities to represent air. The SBP-SAT method represented a cavity by a vacuum and simulates no wave propagation inside present cavities.  

The results indicate that the treatment of cavities, boundaries and interfaces influence the solution. We saw that the Rayleigh surface wave was more exclaimed when the boundary of a cavity was given an explicit treatment with the SBP-SAT method. We also saw that if cavities are represented by air it allows for waves in the computational domain, inside the air, of a much smaller wavelength than that of the waves in the surrounding medium. We argued that these waves of short wavelength can give rise to dispersion errors if they are not resolved appropriately.

The simplification to represent cavities by a vacuum was argued to be valid due to the large difference in density between air and the surrounding rock. This simplification allowed for high order accurate treatment of boundary conditions using a computational grid with grid points aligned to the geometry of the boundary. Future work involves getting a more correct model by extending the SBP-SAT scheme to high order accurate coupling between an elastic material and air.

\end{document}